\documentclass[11pt]{article}
\usepackage{authblk}
\usepackage{geometry}
\usepackage{latexsym}
\usepackage{amsmath}
\usepackage{amsfonts}
\usepackage{amssymb}
\usepackage{graphicx}
\usepackage{color}
\usepackage{parskip}
\usepackage[LGR,T1]{fontenc}
\usepackage[english]{babel}
\usepackage{hyperref}
\usepackage{mathtools}
\usepackage[numbers]{natbib}
\usepackage{tabularx}
\usepackage{subfig}
\hypersetup{colorlinks,breaklinks,linkcolor=blue,urlcolor=red,anchorcolor=red,citecolor=blue}
\usepackage{setspace}
\usepackage{bbm}
\usepackage{appendix}
\setcitestyle{square}

\begin{document}

\title{Skew-Normal Diffusions}

\author[1]{Max-Olivier Hongler\thanks{max.hongler@epfl.ch}}
\author[2]{Daniele Rinaldo\thanks{Corresponding author: d.rinaldo@exeter.ac.uk}}
\affil[1]{\small \'Ecole F\'ed\'erale Polytechnique, Lausanne, Switzerland}
\affil[2]{\small University of Exeter Business School and Land, Environment, Economics and Policy Institute, United Kingdom}

\maketitle

\begin{abstract}$\,$

\noindent We construct a class of  stochastic differential equations driven by White Gaussian noise sources whose solutions can be drawn from skewed Gaussian probability laws, here referred as skew-Normal diffusion (SKN) processes.  The non-Gaussian  character results from implementing a nonlinear and time-inhomogneous drift  constructed  via  ad-hoc changes of probability measure (i.e. Doob's $h$-transform). The SKN processes can be alternatively constructed as dynamic censoring models. While explicitly non-Gaussian, the SKN processes share several properties of Gaussian processes, in particular the invariance under linear transformations. This result allows us to discuss analytically the characteristics of this class of stochastic dynamics. As an illustration, we show how linear noisy monitoring  of SKN processes yields a  solvable  finite dimensional and non-linear stochastic filtering which naturally extends  the  Kalman-Bucy Gaussian case.
\end{abstract}
\noindent\textbf{Keywords:} Stochastic differential equations, nonlinear diffusion processes, skew-Normal distribution, $h$-transform ; dynamic censoring ; nonlinear stochastic filtering

\section{Introduction}

The Skew-Normal probability law $SK(\eta;x)$  extends the Gaussian one by regulating the third moment by means of an extra parameter $\eta$, namely  the {\it skewness}. The general definition of this law is
\begin{equation}
\label{SKNDEFINIT}
\left\{
\begin{array}{l}
SK(\eta;x)dx :=  2\phi(x) \Phi(\eta x)dx, \\\\

\phi(x)dx := \frac{e^{ -  \frac{x^{2}}{2}}}{\sqrt{2\pi }}dx  \quad {\rm and} \quad   {\rm Erf}(x) := \frac{2}{\sqrt{\pi}} \int_{0}^{x} e^{-s^{2}} ds, \\\\
\Phi(\eta x) :=  \frac{1}{\sqrt{2 \pi}}  \left[\int_{- \infty}^{\eta x} e^{- \frac{s^{2}}{2}} ds \right] = \frac{1}{2} \left[1 + {\rm Erf}(\frac{\eta x} {\sqrt{2}} )  \right].

\end{array}
\right.
\end{equation}

\noindent This distribution has been introduced in a seminal paper by Azzalini \cite{azzalini1985class} and subsequently vastly expanded \cite{azzalini2013skew}. One generalisation is known as the extended skew-Normal distribution (ESK)  \cite{azzalini2005skew} and reads:

\begin{equation}
\label{ESKN}
ESK(\eta, x_0;x)dx :=  \phi(x-x_0)\frac{ \Phi(\frac{\eta}{\sqrt{2}} x)}
{ 
\left[  \Phi\left(\frac{\eta x_0}{\sqrt{2(1 + \eta^{2})}} \right)
\right]
}dx.
\end{equation}

\noindent  In Eq.(\ref{ESKN}), the extra $x_0$ parameter offer   flexibility in modelling both moments and tails. From Eqs.(\ref{SKNDEFINIT}) and (\ref{ESKN}), we straightforwardly have $ESK(\eta, 0;x)= SK(\eta, x) $.
The skew-Normal distribution and its extensions find fruitful applications in statistical theory \cite{capitanio2003graphical}, environmental sciences \cite{loperfido2008network}, econometrics \cite{marchenko2012heckman}, time series \cite{jha2023skew}, finance \cite{harvey2010portfolio} and many more. 
Yet,  in the context of   diffusion processes and  Langevin equations, (i.e. stochastic differential equations - SDE - driven by white Gaussian noise), skew-Normal probability laws  remain,  to the best of our knowledge, barely explored. We address this aspect by unveiling a class of Markovian diffusion processes for which the  transition probability densities  (TPD) exhibit skew-Normal laws.  In the sequel, we call this class of processes skew-Normal  diffusions (SKN). \\

The SKN are constructed  via  an ad-hoc change of probability  measure.  This is achieved by  a  logarithmic transformation of the SDE drift which is known as  the $h$-transform. This property  has been studied first  by Doob \cite{doob1957conditional} in the context of the harmonic functions theory \cite{kunita1969absolute,  ito1965transformation, fleming1997asymptotics}. Among the numerous application, this approach plays a  key role to force  diffusion processes to  reach  a given   final state \cite{fleming1977exit, dai1991stochastic, rogers2000diffusions, chetrite2015nonequilibrium}. The SKN processes  in this paper solve SDEs with nonlinear and  time-inhomogeneous drifts. We emphasise that they substantially differ from skew Brownian motions (SBM) which also attracted a sustained attention \cite{ito1963brownian,harrison1981skew}.  The SBM processes  result  from  a reflection at the origin of the standard BM. Specifically, the SBM processes solve SDEs with drifts involving a local time component \cite{ouknine1991skew, lejay2006constructions, appuhamillage2011occupation}, with applications in financial and physical contexts  \cite{decamps2004applications, gairat2017density, bressloff2022probabilistic}.

 \noindent The paper is organised as follows. In Section \ref{s1}, we briefly review $h$-transforms of SDE's. This approach is then used in section \ref{CMEASURE} to present a general result on skewed diffusions and to construct the class of SKN Markov diffusions.   In sections \ref{ALTERNATECON}, we present alternative constructions based on marginal truncations of 2D diffusion processes  and dynamic censoring \ref{DCENSOR} of Brownian motions. In section \ref{PROPO}, we explore some of the exceptional properties enjoyed by the SKN diffusions, in particular their invariance under linear transformations. Finally, in Section \ref{FFILTER}, we show how  the SKN diffusions enable to construct exactly solvable finite dimensional nonlinear filters which naturally extend the well-known Kalman-Bucy case.


\section{H-transforms and skewed diffusive Markov processes}\label{CMEASURE}
\label{htransf}
\label{s1}

Recall briefly some known properties of diffusion processes (see Chapter II of \cite{Borodin}) to be used in the sequel\footnote{For simplicity we  focus on the scalar case: all  results remain valid for $d$-dimensional case}. Consider the  Markov diffusion process in the time interval $[0,T]$ in the filtered probability space $(\mathbb{R},\mathcal{F}, P)$ defined by the SDE \cite{GARDINER}:

\begin{equation}
   dX_t = f(X_t, t) dt +  \sigma(X_t, t) d W_t, \qquad X_0 = x_0. \label{sdeinit}
\end{equation}
where $W_t$ is a standard Brownian motion (BM) \footnote{As  usual $f(x,t)$ to be bounded, continuous and satisfy H\"older conditions in $x$. We further assume $\sigma^2 := \sigma(x,t)^2$ in $C^1(\mathbb{R} \times [0,T])$ and $\sigma^2, \sigma^2_x$ to be bounded and satisfy H\"older conditions in $x$ and $t$. }.  Hence $dW_t$ stands for  the standard White Gaussian noise (WGN) process. As $X_t$  is a Markov process, it is entirely characterised by its transition probability density $ P(x,t| x_0, 0) \in \mathbb{R}^{+}$  which itself 
solves  the Fokker-Planck equation (FPE)

\begin{equation}  \label{kfebase}
\left\{
\begin{array}{l}
\partial_t P(x,t| x_0, 0) =  {\cal L} \left\{  P(x,t| x_0, 0)\right\}, \\\\

{\cal L}\left\{ [\cdot] \right\} := - \partial_x\left\{  f(x,t) [ \cdot] \right\} + \frac{1}{2} \partial_{xx} \left[ \sigma^2(x,t) [\cdot]\right],\\\\
P(x,0| x_0, 0) = \delta( x- x_0)
\end{array}
\right. 
\end{equation}

\noindent with  $\delta(x-x_0)$ being the Dirac probability mass centred at $x_0$. 

\noindent {\bf The} $h$-{\bf transform} \cite{Borodin}.  For  $H(x,t): \mathbb{R} \times [0,T] \rightarrow \mathbb{R}^{+}$ such that:

\begin{equation}  \label{heat1}
\left\{ 
\begin{array}{l}

\partial_t H(x,t) + {\cal L}^{\dagger}  \left\{ H(x,t) \right\}  =0, \\\\ 
{\cal L}^{\dagger} \left\{\cdot \right\} :=   f(x,t)\partial_x [\cdot]  +  \frac{1}{2} \sigma^2(x,t) \partial_{xx} [\cdot].
\end{array}
\right. 
\end{equation}

\noindent    Based on  Eqs.(\ref{sdeinit}) and (\ref{heat1}), consider the process  $\tilde{X}_t$ which  solves the SDE given by

\begin{equation} \label{gensdeh}
d\tilde{X}_t =\big [  f(\tilde{X}_t,  t ) + \sigma^2(\tilde{X}_t,t) \nabla \log H(\tilde{X}_t,t) \big ] dt + \sigma(\tilde{X}_t, t) dW_t ,  \qquad \tilde{X}_0 = x_0.
\end{equation}

\noindent The transformed process  $\tilde{X}_t$ is again Markovian with  TPD $Q(x,t| x_0, 0)$ and the TPDs characterising  $\tilde{X}_t$ and $X_t$  are related as

\begin{equation}\label{qtransf1}
\left\{
\begin{array}{l}
Q(x,t| x_0, 0) = P(x,t| x_0, 0) \frac{ H(x,t)}{H(x_0,0)},\\\\ 
\int_{\mathbb{R}}  Q(x,t| x_0, 0)  dx =  \int_{\mathbb{R}}  P(x,t| x_0, 0) \frac{ H(x,t)}{H(x_0,0)}dx =1.
\end{array}
\right. 
\end{equation}

 \noindent  Eq.(\ref{qtransf1}) shows that  the process $H(X_t,t)/H(x_0,0)$ itself is a positive definite martingale, and the processes  $X_t$ and $\tilde{X}_t$ are related by the following change of probability measure:  
 \begin{equation}
\label{CHPROB}
  dx \mapsto d\mu(x) := \frac{H(X_t,t)}{H(x_0,0)} dx.
\end{equation}

\subsection{A particular class of measure changes}\label{SCMEASURE}

\noindent The previous general framework is now implemented for our purposes with the help of the backward Fokker-Planck equation (BFPE).  To Eq.(\ref{kfebase}), we may associate the BFPE which for  $t\in[0,T]$ reads \cite{GARDINER}:

\begin{equation}
\label{BACKWARD}
\begin{array}{l}
\partial_t   P[x_0,T| x, (T-t)]  + {\cal L}^{\dagger}\left\{ P[x_0,T| x, (T-t)]  \right\}=0.
\end{array}
\end{equation}

\noindent  Define the cumulative probability distribution:

\begin{equation}
\label{HHH}
{\cal H}[x,(T-t)|x_0]:= \int_{-\infty}^{x}   P[x_0,T| \xi, (T-t)] d\xi  = \int_{\mathbb{R}}   P[x_0,T| \xi, (T-t)]   \mathbbm{1} [\xi \leq x ]  d\xi  \in \mathbb{R}^{+},
\end{equation}

\noindent  where $ \mathbbm{1} [\xi \leq x ] $ stands for the Heaviside step  function. By linearity of Eq.(\ref{BACKWARD}),  ${\cal H}[x,(T-t)|x_0]$ itself  solves the BFPE $\partial_t {\cal H}[x,(T-t)|x_0] = {\cal L}^{\dagger}\left\{ {\cal H}[x,(T-t)|x_0]\right\}$. Therefore ${\cal H}[x,(T-t)|x_0]$ satisfies all the necessary conditions to implement a probability change of measures. In this case  Eq.(\ref{qtransf1}) reads:

\begin{equation}
\label{HTRA}
\left\{
\begin{array}{l}
Q(x,t | x_0, 0) = P(x,t | x_0, 0)
\frac{{\cal H}[x,(T-t)|x_0]}{{\cal H}[x_0,T|x_0]}, \qquad t \in [0,T],\\\\
\displaystyle \lim_{t \rightarrow T}Q(x,T | x_0, 0) = P(x,T | x_0, 0) \frac{\mathbbm{1}[x\geq x_0]}{{\cal H}[x_0,T|x_0]} \in [x_0, \infty],
\end{array}
\right. 
\end{equation}
which shows how the $h$-transform constructed in Eq. \eqref{HTRA}  skews the measure $P$ to the right, and at the final time $T$ limits the achievable states of the diffusion to the state space $x \geq x_0$.


\subsubsection{Skew-Normal Diffusions}

\noindent Consider now the BM case, namely $f(x,t)=0$ in Eq.(\ref{sdeinit}), (without loss of generality, we from now on restrict to unit variance):

\begin{equation}
\label{BMM}
\begin{array}{l}
dX_t = dW_t, \qquad X_0=x_0, \quad \Rightarrow \quad p_{{\rm BM}}(x,t |x_0,0) = \frac{e^{- \frac{(x-x_0)^{2}}{2t}} }{\sqrt{2 \pi t}}= \phi\left(\frac{x-x_0}{\sqrt{t} } \right), 
\end{array}
\end{equation}

\noindent and here the solution of the BFPE and its distribution read:

\begin{equation}
\label{BBM}
\left\{
\begin{array}{l}
p_{{\rm BM}}[y,T|x,(T-t)] = \frac{e^{- \frac{(x-y)^{2}}{2(T-t)}}}{\sqrt{2 \pi (T-t)}}, \\\\

{\cal H}[x,(T-t)| y, T ] = \int_{- \infty}^{x}  \frac{e^{- \frac{(\xi-y)^{2}}{2(T-t)}}}{\sqrt{2 \pi (T-t)}} d\xi = \frac{1}{\sqrt{\pi}} \int_{- \infty}^{\frac{x-y}{\sqrt{2(T-t)}}}  e^{z^{2}}dz= \frac{1}{2} \left[ 1 + {\rm Erf} \left( \frac{x-y}{\sqrt{2(T-t)}} \right) \right].
\end{array}
\right.
\end{equation}

\noindent From Eqs.(\ref{BMM}) and (\ref{BBM}), we can establish:

\vspace{0.3cm}
\noindent {\bf Proposition 1.}

\noindent {\it For $t \in [0,T]$, the couple of  dual  Markov diffusion processes solving the SDE:}

\begin{equation}
\label{SDESKEW}
\begin{array}{l}
dX_t = f_{\pm}(X_t, t) dt + dW_t, \qquad X_0=x_0, \qquad t\in [0,T],
\\\\
f_{\pm}(x,t) :=  \partial_x \ln {\cal H}_{\pm}(x,T-t)   =  \partial_x \ln \left[  \Phi\left( \frac{ \pm x}{\sqrt{(T-t)}}\right) \right] ,
\end{array}
\end{equation}

\noindent {\it are  fully  characterised by the dual skew-Normal TPDs $p_{\pm}(x,t|x_0, 0)$ given by}

\begin{equation}
\label{TPDSKEW}
p_{\pm}(x,t|x_0, 0) =
 \phi\left(  \frac{x-x_0}{\sqrt{t}}\right)  \left[ \frac{ \Phi\left( \frac{\pm x}{\sqrt{(T-t)}}\right)}{ \Phi\left( \frac{ \pm x_0}{\sqrt{T}}\right)} \right] =
  ESK\left(\sqrt{\frac{t}{T-t}}, \pm \frac{x_0}{\sqrt{t}};  \pm\frac{x}{\sqrt{t}} \right),  \quad t\in [0,T].
\end{equation}

\vspace{0.2cm}
{\bf Proof of Proposition 1}.

\noindent { The proposition follows directly from the previous results, Eq. \eqref{qtransf1} and the definition Eq.(\ref{ESKN}). Alternatively, the proof can be  obtained by direct substitution of Eq.(\ref{TPDSKEW}) into the FPE:

\begin{equation}
\label{FPSKEW}
\left\{
\begin{array}{l}
\partial_t p_{\pm}(x,t|x_0, 0)  = {\cal L} [ p_{\pm}(x,t|x_0, 0) ], \\\\
{\cal L} [ \cdot] := - \partial_x \left\{   f_{\pm}(x, t)  [\cdot] \right\} + \frac{1}{2} \partial_{xx} [ \cdot].
\end{array}
\right.
\end{equation}

\rightline{$\Box$}

\begin{figure}
    \centering
    \includegraphics[width=11cm]{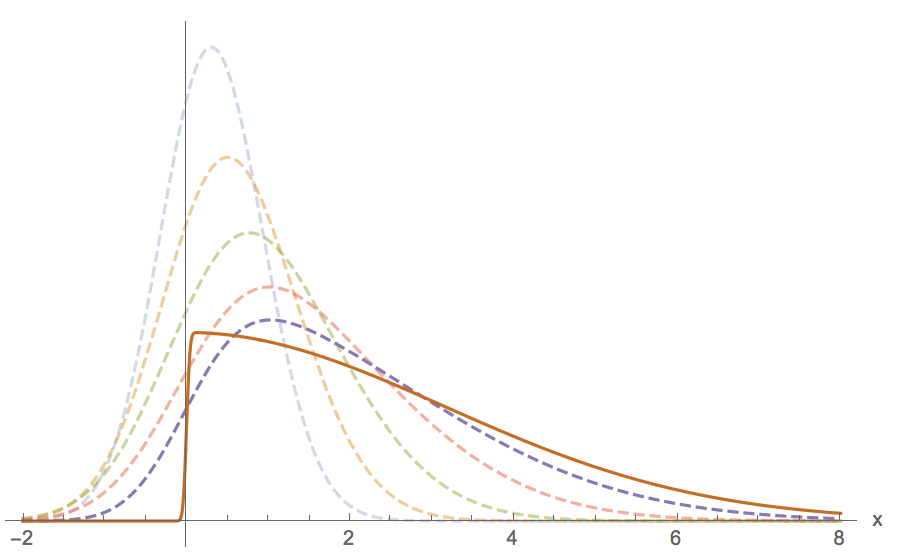}
    \caption{Transient of the right-skewed TPD $p_{\pm}(x,t|0, 0)$ as given by  Eq.(\ref{TPDSKEW}) for the time horizon   $T=10$ and increasing time $t$  from $0.5$ to $T$ (thick red curve). By definition of the   the Heaviside step function   $\mathbbm{1}(x)$, we have  $\lim_{t \rightarrow T} \Phi\left( \frac{\sqrt{2}}{\sqrt{T-t}} \right): = \mathbbm{1}(x)$ thus leading to the half-Normal law (red curve). }
    \label{Fig1}
\end{figure}

  \vspace{0.3cm}
\noindent    {\bf Remark 1}.  For $x_0=0$ implying $\Phi(0)= \frac{1}{2}$,  the  TPD  $p_{\pm}(x,t|0, 0)$  take the  standard SKN  form  Eq.(\ref{SKNDEFINIT}).

\begin{equation}
\label{MATCHO}
p_{\pm}(x,t|0, 0) =  SK\left( \frac{t}{\sqrt{T-t}};  \frac{\pm x}{\sqrt{t}}\right), \qquad t\in [0,T].
\end{equation}

 \noindent  Note that  we have  $p_{{\rm BM}}(x,t |x_0,0) = p_{{\rm BM}}(x-x_0),t |0,0) = \phi\left( \frac{x-x_0}{\sqrt{t}}\right) $.  This translation property  does however not hold for the SKN process:

 \begin{equation}
\label{REMARKIMP}
ESK\left(\sqrt{\frac{t}{T-t}}, \pm \frac{x_0}{\sqrt{t}};  \pm\frac{x}{\sqrt{t}} \right)  \neq SK\left( \frac{t}{\sqrt{T-t}};  \frac{\pm (x-x_0)}{\sqrt{t}}\right).
\end{equation}

\noindent From Eq.(\ref{SDESKEW}), we  note  that $\left[ {\cal H}_{-}(x,T-t) + {\cal H}_{+}(x,T-t)\right] =1$. 

\vspace{0.3cm}
\noindent {\bf Remark 2}. The TPD Eq.(\ref{TPDSKEW})  is defined only  for finite  time horizon $t \in [0,T]$. One may  wonder whether one could construct a SKN process for unrestricted time horizon  $t\in \mathbb{R}^+$.  In Appendix B, this question  is briefly addressed using the  as $h$-transform involving imaginary error function ${\rm Erfi}(x)$.  There one shows that such  process exists  only on positive state space.


\section{Alternative constructions of skew-Normal diffusions} \label{ALTERNATECON}

\noindent Section \ref{CMEASURE} shows that  SKN diffusion processes  result directly from a change of probability of the BM. Here we present a set of  alternative interpretations.
 \subsection{The SKN as a truncated  marginal of bi-dimensional  Gaussian process}
\label{skn}

\noindent The  centred binormal probability law is defined as:

 \begin{equation}
 \label{BINGAU}
 \begin{array}{l}
 {\cal N} (x,y) dx \, dy= \frac{1}{2\pi \sigma_1 \sigma_2 \sqrt{1- \rho^{2}} }  \exp \left\{-\frac{(x/\sigma_1)^{2}  + 2 \rho (x/\sigma_1) (y/\sigma_2) - (y/\sigma_2)^{2}}{2 (1- \rho^{2} )} \right\} dx\, dy ,
 
 \end{array}
 \end{equation}
 \noindent with variances  $\sigma_1^2$ and $\sigma_2^2$ and correlations $\rho \in [-1,+1]$. Define a couple of dual  truncated marginals as: 
 $\Pi^{\pm}(x)$ as :

 \begin{equation}
\label{TRUNCATED}
 \left\{ 
 \begin{array}{l}
 \Pi^{\pm} (x)dx:=\left[ \int_{\mathbb{R}^{\pm}}  {\cal N} (x,y) dy \right] dx = \frac{e^{-(x/ \sqrt{2}\sigma_1)^{2}}}{\sqrt{2\pi  \sigma_1^{2}
 } }   d\mu_{\pm}(x) =  \phi(\frac{x}{\sigma_1})  d\mu_{\pm}(x)
 \\\\

 d\mu_{\pm}(x) :=\frac{1}{2} \left[ 1 \pm  {\rm Erf} \left( \frac{\rho x}{\sqrt{2 \sigma_1^{2} (1- \rho^{2}) }} \right)  \right] dx, 
 \end{array}
 \right.
 \end{equation}

 \noindent  with  $\Pi^{-}(x) + \Pi^{+}(x) = \phi(\frac{x}{\sigma_1})$. Focusing on  SKN diffusion processes, let us  introduce  a natural dynamic extension  of Eqs.(\ref{BINGAU})  and (\ref{TRUNCATED}) by writing  a bi-dimensional diffusion process  $(X_t,Y_t) \in \mathbb{R}^2 $, $t \in \mathbb{R}^{+}$  and 
\begin{equation}
\label{DYNAMO}
X_t = W_{1,t} \quad {\rm and} \quad  Y_t=  \sqrt{a(t)} W_{1,t} + \sqrt{b(t)} W_{2,t},  \qquad a(t) \geq 0\,\, {\rm and} \,\,  b(t) \geq 0,
\end{equation}

\noindent where $W_{1,t}$ and $W_{2,t}$ are independent standard BM's and $a(t), b(t)$  time-dependent functions. By definition of the BM, we have:

\begin{equation}
\label{NDEF}
\sigma^{2}_1(t) = t \quad \sigma^{2}_2(t) =t[a(t) + b(t)] \quad {\rm and} \quad \rho^{2}(t) = a(t)t.
\end{equation}

\noindent Observe at this point that $a(t)$ cannot be chosen arbitrarily since  we must have   $\rho^{2}(t) \in [-1,+1]$. From  Eqs.(\ref{BINGAU}) and (\ref{NDEF}),  the resulting  stochastic  process $(X_t, Y_t)$ is a  Gaussian Markov diffusion on $\mathbb{R}^{2}$ with TOD:

 \begin{equation}
 \label{BIN}
 \begin{array}{l}
p(x,y,t  |x_0, y_0,0) = \frac{1}{2\pi t \sigma_2 (t)\sqrt{1- \rho^{2}(t)} }  \exp \left\{-\frac{[(x-x_0)^{2}/2t]  + 2 \rho(t) [x-x_0/\sqrt{t}] [(y-y_0)/\sigma_2(t))] -[ (y-y_0)^{2}/ \sigma^{2}_2(t))]}{2 (1- \rho^{2} (t))} \right\} dx\, dy ,
 \end{array}
 \end{equation}

\noindent Fix now   $x_0= y_0=0$ and   as  in Eq.(\ref{TRUNCATED}), define:

\begin{equation}
\label{PIT}
\left\{
\begin{array}{l}
 \int_{\mathbb{R}} p(x,y,t  |0, 0,0) dy := \Pi^{\pm}(x,t) = \frac{ e^{-\frac{x^{2}}{2t}} }{2\sqrt{2\pi t}}  
\left\{ 1 \pm {\rm Erf}\left[ Z(t) x\right] \right\}  dx=\frac{1}{2}  SK\left( \pm \sqrt{2 t}Z(t) ;\frac{x}{\sqrt{t}} \right)dx , \\\\
Z(t) :=  \frac{\rho(t) \,  }{\sqrt{2 t[1-\rho^{2}(t) ] }} \quad {\rm with} \quad  \rho(t) \in [-1,+1].
\end{array}
\right.
\end{equation}

\noindent  Let us now study the time evolution of  $\Pi^{\pm}(x,t)$.

\vspace{0.3cm}
\noindent{\bf Proposition 2.} {\it Given the   correlations $\rho(t) \in [-1, +1]$, the SKN law $ \Pi^{\pm}(x,t) $ Eq.(\ref{PIT}) solves  the  parabolic pde:}

\begin{equation}
\label{MATCH1}
\left\{
\begin{array}{l}
\partial_t\Pi^{\pm}(x,t) = - \partial_x \left[f^{\pm}_{\psi}(x,t) \Pi^{\pm}(x,t) \right]  +   \frac{1}{2} \partial_{xx} [\Pi^{\pm}(x,t)], \\\\
f^{\pm}_{\psi}(x,t)   :=  \psi(t) \partial_x  \ln \left\{  1 \pm {\rm Erf}[Z(t)x] \right\}.
 \end{array}
\right.
\end{equation}

\noindent {\it where $\psi(t) $ is given by:}

\begin{equation}
\label{RHOPROP}
\begin{array}{l}

 \psi(t) = \frac{1}{2} + t \partial_t \left[ \ln \rho(t) \right] .
 
 \end{array}
\end{equation}

\vspace{0.3cm}
\noindent{\bf Proof of Proposition 2:}  By direct verification, as detailed  in Appendix A.

\rightline{$\Box$}

\vspace{0.3cm}

\noindent{\bf Remark 3}.

\noindent The density  $\Pi^{\pm}(x,t) = \int_{\mathbb{R}} p(x,y,t  |0, 0,0) dy$ is a projection (truncation mechanism)  of a  bi-normal Markov diffusion for the special initial condition $x_0 = y_0 =0$. Such a projection mechanism destroys information and so the Markov property of  the nominal  2D diffusion process is  generally not be preserved for arbitrary $\psi(t)$. In Corollary 2 below, it is shown that  the  information  loss can be exactly  balanced  by an  information  gain obtained by an ad-hoc correlation  increase  between the components of the 2D diffusion Eq.(\ref{BIN}).  

\noindent{\bf Corollary 2}

\noindent {\it  For finite time horizon  $t\in [0,T]$ and correlations  $\rho^{2}(t) = t/T$, Eq.(\ref{MATCH1}) and  Eq.(\ref{SDESKEW})  are identical. For arbitrary  $x_0 \in \mathbb{R}$, both are solved by $ ESK\left(\sqrt{\frac{t}{T-t}}, \pm \frac{x_0}{\sqrt{t}};  \pm\frac{x}{\sqrt{t}} \right) $ which is the TPD characterising  a Markov diffusion process.}

\vspace{0.3cm}
\noindent{\bf Proof of Corollary 2:}  { For correlations $\rho^{2}(t) = t/T$, Eq.(\ref{RHOPROP}) implies that $\psi(t) =1$. Accordingly the drift component  in the pde  Eq.(\ref{MATCH1}) coincides with the one in  Eq.(\ref{SDESKEW}) of  Proposition 1.  Accordingly, the solution of the corresponding FPE Eq.(\ref{FPSKEW}) is given by $ ESK\left(\sqrt{\frac{t}{T-t}}, \pm \frac{x_0}{\sqrt{t}};  \pm\frac{x}{\sqrt{t}} \right) $ which fully characterises a Markov diffusion process. }
 
\rightline{$\Box$}

Note that setting $\rho_t = \rho \in [0,1)$ (constant correlations, implying $a(t) = \rho, b(t) = \sqrt{1 - \rho^2}, \psi(t) = 1/2$),
yields a skew-Normal diffusion $X_t^\pm$ with skewness proportional to $\pm 1/\sqrt{t}$. This specific time-dependence is a suitable feature for stochastic volatility models in finance: for example, multiple empirical studies of the equity skew \cite{carr2003finite, azzone2022additive} observe how the S\&P 500 short-time skew exhibits the same property as the SKN.


\subsection{SKN and dynamic censoring}\label{DCENSOR}
\label{s5.3}

\noindent  Assuming $x$ as  a centered Gaussian random variable with unit variance and law $\phi(x)$ given in   Eq.(\ref{SKNDEFINIT}),  one  defines a couple of  one-sided truncated (also referred as censored) probability laws $T_{[a]}(x)$ and $T^{[b]}(x)$ as:

\begin{equation}
\label{TRUNCA}
\left\{
\begin{array}{l}
T_{[a_{\downarrow}]}(x): = {\rm Prob}\left\{x | x> a_{\downarrow}\right\}  =  \int_{a_{\downarrow}}^{ \infty} \phi(x) dx  =  +  \frac{2\phi(x) }{1 -{\rm Erf}\left( \frac{a_{\downarrow}}{\sqrt{2}} \right)   }, 
\\\\ 
T_{[a_{\uparrow}]}(x): = {\rm Prob}\left\{x | x< a_{\uparrow} \right\}  =  \,  \int_{- \infty}^{a_\uparrow} \phi(x)dx =  -\frac{2\phi(x) }{1 -  {\rm Erf}\left( \frac{a_{\uparrow}}{\sqrt{2}} \right)   }, 
\end{array}
\right.
\end{equation}

\noindent  Writing $a_{\downarrow}= -(x+  \epsilon)$ and   $ a_{\uparrow}= + (x+ \epsilon)$ in Eq.(\ref{TRUNCA}) and recalling that  $ {\rm Erf}\left[- \frac{x}{\sqrt{2}} \right]= - {\rm Erf}\left[ \frac{x}{\sqrt{2}}\right]$, we have:

\begin{equation}
\label{FORMALITY}
\left\{
\begin{array}{l}
\displaystyle \lim_{\epsilon \rightarrow 0} T_{[-(x + \epsilon)]}(x) =  +  \lim_{\epsilon \rightarrow 0}\left[ \frac{2\phi(x) }{1 - {\rm Erf}\left(- \frac{x+ \epsilon }{\sqrt{2}} \right)   }\right] =  \sqrt{2} \partial_x \ln  \left[1 + {\rm Erf}\left( \frac{x}{\sqrt{2}} \right)\right],

 \\\\
 \displaystyle \lim_{\epsilon \rightarrow 0} T_{[+(x+ \epsilon)] }(x) =  -  \lim_{\epsilon \rightarrow 0}  \left[ \frac{2\phi(x) }{1 - {\rm Erf}\left(+ \frac{x+ \epsilon}{\sqrt{2}} \right)   }\right] = \sqrt{2} \partial_x \ln  \left[1- {\rm Erf}\left( \frac{x}{\sqrt{2}} \right)\right].
\end{array}
\right.
\end{equation}

\noindent  
 The manifest  similarity of Eq.(\ref{FORMALITY}) with the drifts terms  in Eq.(\ref{MATCH1})  suggests that SKN processes can be understood as an   ad-hoc  dynamic truncation  of a Brownian motion. Indeed, by drawing the increments of a BM from a truncated Gaussian law we can directly generate the SKN drifts. By choosing  $a_{\uparrow, (\downarrow) } =  \mp  Z(t) X_t$ in Eq.(\ref{FORMALITY}) we directly  generate the drifts appearing in Eq.(\ref{MATCH1}). Summarising we  have:

\textbf{Proposition 3: } {\it The SKN biases induced by the  dual drifts in Eq.(\ref{MATCH1}) are  generated by the dynamic censoring of a BM}:
 \begin{equation}
\label{FORMALTRUNCA}
\left\{
\begin{array}{l}
f_{-}(X_t,t) = 
    \psi(t)  T_{\left[ -Z(t) X_t \right]}(X_t), \\\\
   
   f_{+}(X_t,t) = 
    \psi(t)  T_{\left[ + Z(t) X_t \right]}(X_t) .
\end{array}
\right.
\end{equation}

\vspace{0.3cm}

\textbf{Proof of Proposition 3 }.

\noindent The assertion follows directly from Eqs.(\ref{MATCH1}), (\ref{TRUNCA}) and (\ref{FORMALITY}).

\rightline{$\Box$}

\noindent {\bf Remark 4}:

\noindent The characteristic  $\frac{\phi(x)}{\Phi(\pm x)}$ structure of the  drifts Eq.\eqref{MATCH1} is typical  in reliability theory known as the {\it hazard rate function}  and in statistics known as the  {\it inverse Mills ratio}. It naturally arises in the context  of  the censoring of random variables well known in econometrics  \cite{hayashi2011econometrics}. In particular,  the relevance of the skew-Normal laws in a static censoring mechanism is  well known  \cite{heckman1979sample, azzalini2019sample}. 
Here we have instead a dynamic  extension of this mechanism. 
For all $t \in [0,T]$, the skewness $Z(t)$ increasingly expands the censored part of the domain of $X_t$, ultimately leading to a censoring on the half-line $X^+_{T} \geq x_0$ (respectively $X_T^- \leq x_0$).

\section{Properties of the SKN diffusion process}\label{PROPO}

\noindent Since the  SKN diffusion processes can be derived from a Brownian motion via a change of probability measure, they retain several remarkable properties of the Brownian motion itself.

\subsection{Moments of the SKN diffusion process} \label{MOMO}

\noindent  \noindent {\bf Corollary 3}

\noindent {\it The generating function for  the moments of the SKN process Eq.(\ref{TPDSKEW})  reads:}

\begin{equation}
\label{MOMENTSG}
G_{\pm}(k,t) :=    \int_{\mathbb{R}} e^{k(x-x_0)}    ESK\left(\sqrt{\frac{t}{T-t}}, \pm \frac{x_0}{\sqrt{t}};  \pm\frac{x}{\sqrt{t}} \right)  dx =  \frac{\left[ e^{\frac{k^{2} t}{2}}\pm {\rm Erf}\left( \frac{ x_0 -kt}{\sqrt{2T}}\right)\right] }{\left[1\pm {\rm Erf}\left( \frac{ x_0}{\sqrt{2T}}\right)\right]}.
\end{equation}

\vspace{0.3cm}
\noindent {\bf Proof of Corollary 3}

\noindent By direct integration, see Appendix C.

 \rightline{$\Box$}

 \noindent The first two moments follow easily by using Eq.(\ref{MOMENTSG}):
 
 \begin{equation}
\label{ONETWO}
 \left\{
 \begin{array}{l}
 \mathbb{E} \left\{(x-x_0)  \right\}  = \frac{d}{dk} G_{\pm}(k,t) \mid_{k=0} =t f_{\mp}(x_0), \\\\
 
  \mathbb{E} \left\{(x-x_0)^{2}  \right\}  = \frac{d^{2}}{dk^{2}} G_{\pm}(k,t) \mid_{k=0} = t \left[1 - \frac{t}{T}  x_0 f_{\mp}(x_0) \right], 

 \end{array}
 \right.
\end{equation}
the latter result showing how SKN processes behave diffusively for $x_0 = 0$ (the variance of the process increases linearly in time), and sub-diffusively \cite{alves2016characterization, massignan2014nonergodic} as $x_0$ increases (for $f_+(x_0)$,  viceversa for $f_-(x_0))$. 
\subsection{Linear transformations of the ESKN processes} \label{INVARO}

\noindent  A fundamental property of the Gaussian  law is its invariance under linear transformations. This exceptional property is preserved in SKN processes. Let us show this by considering the degenerate 2-D diffusion process n  $\mathbb{R}^{2}$:

 \begin{equation}
\label{DRIVE2}
\left\{
\begin{array}{l}
dX_t =  - \lambda X_t  dt+ d U_t, \qquad X_0 =x_0,\\\\ 

dU_t =  \partial_u  \ln\left\{ 1 +  {\rm Erf}\left[ \frac{u}{\sqrt{2(T-t)}}\right]  \right\}_{u= U_t}dt+ dW_t, \, \,\, U_0 =0 \, \,\,  {\rm for} \,\,\, t\in [0,T].
\end{array}
\right.
\end{equation}
  
\noindent  Eq.(\ref{DRIVE}) is a linear transformation of the Markov  SKN process $U_t$. Physically speaking, the process $X_t$  Eq.(\ref{DRIVE}) can be seen as  a  (linear)  Ornstein-Uhlenbeck process driven by a SKN noise source $dU_t$. The  bi-variate Markov process $(X_t, U_t)$ is thus completely characterised by its TPD  $p(x, u,t |x_0, 0)$ and we have:

\noindent {\bf Corollary 4}

\noindent {\it The $X_t$-marginal density of the 2D degenerate diffusion process Eq.(\ref{DRIVE})  belongs to the  ESKN family and reads:}

\begin{equation}
\label{MARGO1}
P_M (x,t|x_0, 0) := \int_{\mathbb{R}} p(x, u,t |x_0, 0)  du= \frac{e^{- \frac{\left(x-x_0e^{-\lambda t}\right)^{2} } {2 r(t)}} }{\sqrt{2 \pi r(t)}}
 \left\{
 1 + {\rm Erf} \left[ \frac{ s(t)\left[ x-x_0e^{-\lambda t}\right]}
{
\sqrt{2 \left[   r^{2}(t)T - r(t)s^{2}(t) e^{-4\lambda t} \right]}
}
 \right] 
\right\}, 
\end{equation}
\noindent{\it with  $r(t) = \frac{1-e^{-2\lambda t}}{2 \lambda} \quad {\rm and} \quad s(t) = \frac{1-e^{-\lambda t}}{\lambda}$.}

\vspace{0.3cm}
\noindent{\bf Proof of Corollary 4}

\noindent The full calculation is displayed in  Appendix D.

\rightline{$\Box$}

\noindent From Eq.(\ref{MARGO1}),  we straightforwardly have  $\displaystyle \lim_{\lambda \rightarrow 0} r(t) =   \lim_{\lambda \rightarrow 0} s(t) = t$. So it is immediate to see that  in the $\lambda =0$ limit, Eq.(\ref{MARGO1}), reduces to:

\begin{equation}
\label{MARGOLIM}
\displaystyle \lim_{\lambda \rightarrow 0}  P_M (x,t|x_0, 0) = \frac{e^{- \frac{ (x-x_0)^{2} } {2t}} }{\sqrt{2 \pi t}} \left\{ 1 + {\rm Erf} \left[ \frac{x-x_0}{\sqrt{2 (T-t)}} \right] \right\}. 
\end{equation}


\subsection{The BM as the superposition of  SKN diffusions}\label{LUMP}

\noindent From  Proposition 1 and in particular Eq.(\ref{TPDSKEW}), we can write:

\begin{equation}
\label{SUPERPO}
\Phi\left[\frac{-x_0}{\sqrt{T}}  \right]   p_{-} (x,t |x_0) + \Phi\left[\frac{+x_0}{\sqrt{T}} \right]  p_{-} (x,t |x_0)= \frac{e^{- \frac{(x-x_0)^{2}}{2t} } }{\sqrt{2\pi t}}. 
\end{equation}

\noindent  Intuitively, Eq.(\ref{SUPERPO}) indicates that the BM's TPD with  initial position $x_0$ can be represented as a weighted superposition of left- and right skewed  extended SKN processes with opposite chiralities. 
This is an explicit illustration of Markov functions as introduced in the seminal work of   L. C. Rogers and  J. Pitman  \cite{rogers1981markov}. Directly inspired by their Example 2, define a  state space ${\cal S} := \left\{-1, +1 \right\}\times \mathbb{R}$ on which evolves a continuous time Markov  stochastic process: $\Xi_t = ({\cal B}, X_t)$ characterised   as follows:

\begin{itemize}
\item[]i) \,
Consider a Bernoulli process ${\cal B}_{t}$ with realisations $\left\{-1,+1 \right\}$. At  initial time $t=0$,  ${\cal B}_{0}$ delivers either  the  outcome  $b_0=  \left\{ -1, + 1\right\}$ with respective probabilites $\Phi\left[\frac{\pm x_0}{\sqrt{T}}  \right] \in [0,1]$.

  \item[]ii)  For  $t>0$,  the diffusive  component of the process  $\Xi_t$ is given by :

\begin{equation}
\label{DRIFT}
\left\{
\begin{array}{l}
dX_t = [ f_{b_0}(X_t,t)]dt + dW_t,  \qquad X_0=x_0,\\\\

f_{b_0}(x,t)  =  \partial_x\ln  \left\{ 1 + b_0   {\rm Erf} \left[  \frac{x}{\sqrt{2(T-t)}}\right] \right\}, \qquad b_0\in  \left\{   -1 +1  \right\} , \\\\

\partial_t p_{b_0}(x_0,0 |x,t) = - \partial_x \left[ f_{b_0}(x,t)  p_{b_0}(x_0,0 |x,t) \right] + \frac{1}{2} \partial_{xx} p_{b_0}(x_0,0 |x,t) .
\end{array}
 \right.
\end{equation}

\noindent where $f_{b_0}(x,t)$ and  $p_{b_0}(x_0,0 |x,t) =p_{\pm}(x,t |x_0,0)$  are the drifts respectively the TPDs given in  Eq.(\ref{TPDSKEW}). This shows that  the scalar BM  initiated at $x_0$   results  itself from the superposition of  a couple of SKN's  with opposite biases. In other words, the BM itself is a  lumped version on of  the Markov process $\Xi_t $ itself evolving on the enlarged probability  state  ${\cal S}$.
\end{itemize}

\section{Illustration - finite dimensional  nonlinear   stochastic filtering}\label{FFILTER}

\noindent In section \ref{INVARO} we saw that the SKN are invariant under linear transformations. This offers the exceptional opportunity to derive an  exact and explicit solution for a non-linear filtering problem.  

\noindent Let us briefly recall first the  basic problematic  of stochastic filtering. The position $X_t$ of a  stochastic process  is monitored by  a measurement device itself corrupted by WGN. The objective of stochastic filtering is to construct a real time and optimal  (in the mean square sense) estimate $\hat{X}_t$ of the actual  position $X_t$. The utmost  engineering relevance of this general  problematic still is the source of an abundant  literature. Basics references  of direct use  in the sequel are \cite{OKSENDAL, SCHUSS}. The special situation   where  the system's dynamic and the monitoring process both  are  linear is known as the Kalman-Bucy filtering  problem (KFP). The KFP can be solved explicitly, since linearity preserves the  Gaussian character of the driving  noise sources. Accordingly, for the KFP  a finite number of estimate moments of  are sufficient to fully solve the filtering problem.  The SKN  invariance under linear transformations  discussed in section \ref{INVARO},  suggests that observing a SKN process via a linear monitoring platform also   offers the  possibility to  solve exactly  a   3-parameters non-Gaussian filtering problem ; this is discussed in this section.  Start with the general filtering problem (FIP):

\begin{equation}
\label{FILTER1}
\left\{
\begin{array}{l}
dX_t = f(X_t,t) dt+  dW_{1,t}, \quad\,  X_0 =x_0,\qquad {\rm (system \,\, evolution)}, \\\\

d{\cal Z} _t = \omega X_t dt+ \sigma dW_{2,t}, \qquad {\cal Z}_0 = x_0,\,\, \qquad {\rm (monitoring \,\, process)},

\end{array}
\right.
\end{equation}

\noindent where $dW_{1,t} $ and $dW_{2,t} $ stand for independent WGN's and $\sigma$ is a positive constant noise amplitude.  Eq.(\ref{FILTER1}) synthesises  a $\sigma dW_{2,t}$ noisy  linear observation $\omega X_t,  \,\, \omega > 0$  of the system $X_t$\footnote{Without loss of generality the driving noise $dW_{1,t}$-variance  of the system  is chosen to be unity.}.The filtering problem (FIP) consists  in constructing a probability measure  $\varphi(x, t | {\cal Z}_{\tau})$ to estimate   $X_t$. As it is clear from  the notation, the  estimate will be  based on the full observation history $ {\cal Z}_{\tau}$ for $\tau\in [0, t]$. Accordingly, $\varphi(x, t |  {\cal Z}_{\tau})$ is itself  a  random probability  measure which characterises $\hat{X}_t$ and it is  steadily adapted in real time. Numerous distinct  estimators are obviously possible  and here one  focuses on the optimal  one. Hence $\varphi(x, t |  {\cal Z}_{\tau})$ which  minimises  the  quadric errors:

\begin{equation}
\label{ERROR}
 |\hat{e} |^{2}(t):
 = \int_{\mathbb{R}} | x- \hat{x} |^{2} \varphi(x, t |  {\cal Z}_{\tau}) dx,
\end{equation}

\noindent  For  the FIP defined by  Eq.(\ref{FILTER1}), general stochastic filtering  theory establishes that  $\varphi(x, t |  {\cal Z}_{\tau})$  obeys  the {\it  Kushner-Stratonovich} (KS) stochastic pde: 
 
 \begin{equation}
\label{K}
\left\{
\begin{array}{l}
d \varphi(x, t |  {\cal Z}_{\tau}) = {\cal L} [\varphi(x, t |  {\cal Z}_{\tau}) ] dt + \frac{\omega}{\sigma^{2}} \left\{ x-\mathbb{E} \left[ \hat{X}(t) \right] \right\} \left\{   d {\cal Z}_t -\mathbb{E} \left[ \hat{X}(t) \right] )\right\}  \varphi(x, t |  {\cal Z}_{\tau}), \\\\

 {\cal L} [\cdot] = -\partial_x\left[ f(x) [\cdot] \right]  + \frac{1}{2} \partial_{xx} [\cdot],\\\\

\mathbb{E} \left[ \hat{X}(t) \right] := \int_{\mathbb{R}}  x \varphi(x, t |  {\cal Z}_{\tau})  dx.
 \end{array}
 \right.
\end{equation}
\noindent The   KS stochastic  pde is intrinsically non-linear since  it involves  the real time expectation  $\hat{x}(t)$.  It is remarkable that the un-normalised function $\phi(x,t|Z_t)$ defined as:

\begin{equation}
\label{PHI}
\varphi(x, t |  {\cal Z}_{\tau})  := \frac{\psi(x,t| {\cal Z}_t)}{\int_{\mathbb{R}} \psi(x,t| {\cal Z}_t) dx}. 
\end{equation}

\noindent follows a linear  evolution known as  the Zakai equation which  reads \cite{SCHUSS}:

\begin{equation}
\label{Z}
\left\{
\begin{array}{l}
d \psi(x, t |  {\cal Z}_{\tau}) = {\cal L} [\psi(x, t |  {\cal Z}_{\tau}) ] dt + \frac{\omega x}{\sigma^{2}}[\psi(x, t |  {\cal Z}_{\tau}) ] d {\cal Z}_t \qquad \qquad  \qquad ({\rm Statonovich\,\, interpretation}),  \\\\ 
 d \psi(x, t |  {\cal Z}_{\tau}) =  {\cal L}  [\psi(x, t |  {\cal Z}_{\tau})  ] dt - \frac{\omega^{2} x^{2}}{2 \sigma^{2}}  [\psi(x, t |  {\cal Z}_{\tau})  ] dt + \frac{\omega x}{\sigma^{2}}[\psi(x, t |  {\cal Z}_{\tau}) ] d {\cal Z}_t \quad ({\rm Ito\,\, interpretation}).
\end{array}
\right.
\end{equation}

 \noindent Apply these FIP concepts for systems $X_t$ and $\tilde{X}_t$ related via an   $h$-transformation.  In relation with Eq.(\ref{FILTER1}), consider the  FIP with dynamics:

\begin{equation}
\label{FDOOB2}
\left\{
\begin{array}{l}
d\hat{X}_t= \left[ f(\hat{X}_t,t) + \left\{ \partial_{x} \ln [ H(x,t)]\right\}_{x= \hat{X_t}} \right] dt + dW_{1,t}, \quad \hat{X}_0= x_0,
\\\\
d {\cal Z}_t = \omega \hat{X}_t dt+ \sigma dW_{2,t}, \qquad  {\cal Z}_0 = x_0
\end{array}
\right.
\end{equation}
where $\partial_t  H(x,t) + {\cal L} [ H(x,t)]=0$ as introduced in section \ref{htransf}. The modified evolution  Eq.(\ref{FDOOB2})  leads to the modified Zakai dynamics:

\begin{equation}
\label{Z2}
\left\{
\begin{array}{l}
d\hat{\psi}(x, t |  {\cal Z}_{\tau})  =  \hat{{\cal L}}[ \hat{\psi}(x, t |  {\cal Z}_{\tau}) ]dt +\frac{\omega x}{\sigma^{2}}[\hat{\psi}(x, t |  {\cal Z}_{\tau}) ] d {\cal Z}_t, \qquad ({\rm Statonovich\,\, interpretation}), \\\\
\hat{{\cal L}}[\cdot] = - \partial_x\left\{ [f(x,t) + \partial_x \ln [H(x,t)] [\cdot] \right\} + \frac{1}{2} \partial^{2}_{xx} [\cdot]
\end{array}
\right.
\end{equation}

\subsection{SKN and finite-dimensional nonlinear filters}

\noindent The Kalman-Bucy FIP results  when $f(x,t)$ is linear in $x$. 
In this case the FIP involves  a 2D Gaussian diffusion process and $\varphi(x, t |  {\cal Z}_{\tau}) $  is itself Gaussian and its moments can be explicitly calculated. This is a finite dimensional FIP since  the first two moments fully characterised the filter. This situation can be generalised for certain classes of dynamics, an illustration is now exhibited.

\noindent {\bf Proposition 4}. {\it For the BM (i.e.  $f(x,t) =0$), the  un-normalised Zakai filters $\psi(x, t |  {\cal Z}_{\tau})$ and $\hat{\psi}(x, t |  {\cal Z}_{\tau})$ in Eqs.(\ref{Z}) and (\ref{Z2}) are related as:}

\begin{equation}
\label{CMZ}
  \hat{\psi}(x, t |  {\cal Z}_{\tau}) = H(x,t) \psi(x, t |  {\cal Z}_{\tau}) .
\end{equation}

\vspace{0.3cm}
\noindent {\bf Proof of Proposition 4:}

Take $f(x,t) =0$ implying ${ \cal L }[\cdot]  = \frac{1}{2} \partial_{xx}[\cdot]$  with  initial condition $x_0=0$ and $\omega = \sigma =1$ in Eqs.(\ref{Z}) and (\ref{FDOOB2}) and we obtain:

\begin{equation}
\label{DEMO1}
\left\{
\begin{array}{l}
\partial_t H(x,t) +  \frac{1}{2} \partial_{xx} H(x,t) =0, \\\\
d\hat{\psi}(x, t |  {\cal Z}_{\tau})  =- \partial_x \left\{ \left[  \partial_x \ln H(x,t) \right] \hat{\psi}(x, t |  {\cal Z}_{\tau})  \right\}dt  + \frac{1}{2} \partial_{xx}[ \hat{\psi}(x, t |  {\cal Z}_{\tau}) ]dt +\frac{\omega x}{\sigma^{2}}[\hat{\psi}(x, t |  {\cal Z}_{\tau}) ] d {\cal Z}_t, 
\end{array}
\right.
\end{equation}

\noindent  Omitting the arguments,  write  $\hat{\phi} = H \phi$ verify that it solves Eq.(\ref{DEMO1}), namely:

\begin{equation}
\label{DEMO2}
 \psi \left(d H +  \frac{1}{2}\partial_{xx} H\right) - H \left(d\psi  - \frac{1}{2} \partial_{xx} \psi -  \frac{\omega x}{\sigma^{2}}\psi d {\cal Z}_t,  \right) =0.
 \end{equation}
 
 \rightline{$\Box$}

\noindent  When $f(x,t)=0$,  Zakai's Eq.(\ref{Z}) is solved by  a Gaussian, (see Example 6.2.10 in \cite{OKSENDAL} and  section 3.6 in \cite{SCHUSS}) and we have:

\begin{equation}
\label{KFB1}
\left\{
\begin{array}{l}
\psi (x, t |  {\cal Z}_{\tau})=  e^{-\frac{(x- \hat{X}_t)^{2}}{2 S(t)} } \quad {\rm and} \quad  \varphi (x, t |  {\cal Z}_{\tau})= \frac{e^{-\frac{(x- \hat{X}_t)^{2}}{2 S(t)} } 
}{\sqrt{2 \pi S(t) }} 
 \\\\ 
S(t) = \tanh(t) \quad {\rm and } \quad  \hat{X}_t  = \frac{
\int_{0}^{t} \sinh(s) d {\cal Z}_t}
{\cosh(t)}.

\end{array}
\right.
\end{equation}

\noindent In Eq.(\ref{CMZ}), identify $H(x,t) ={\cal H}(x,T-t) =  \Phi\left( \frac{\pm x}{\sqrt{(T-t)}}\right)$ as given by Eq.(\ref{SDESKEW}) and the resulting  FIP  Eq.(\ref{FDOOB2}) is finite dimensional. Hence  using Eq.(\ref{CMZ}), we have immediately a couple of FIPs:

\begin{equation}
\label{FILTEREX}
\left\{
\begin{array}{l}
\hat{\psi} _{\pm}(x, t |  {\cal Z}_{\tau})= \psi (x, t |  {\cal Z}_{\tau})\Phi\left( \frac{\pm x}{\sqrt{(T-t)}}\right),  \qquad \qquad t\in [0,T], \\\\
\hat{\varphi}_{\pm} (x, t |  {\cal Z}_{\tau}) =  \varphi (x, t |  {\cal Z}_{\tau}) \frac{\Phi\left( \frac{\pm x}{\sqrt{(T-t)}}\right)
}{
\Phi\left( \frac{\pm \hat{X}_t \sqrt{\tanh(t)}}{\sqrt{\tanh(t) + (T-t)}}\right)},

\end{array}
\right.
\end{equation}

\noindent where the normalisation factor of   $\hat{\varphi}_{\pm} (x, t |  {\cal Z}_{\tau}) $ follows by using  entry 2.7.1-6 of  \cite{korotkov2020integrals}. Eq.(\ref{FILTEREX}) is an  illustration belonging to the class of  finite  dimensional filters first discussed by V. E.  Benes \cite{BENES}.

\bibliographystyle{unsrtnat}
\bibliography{skn}

@article{dai1991stochastic,
  title={A stochastic control approach to reciprocal diffusion processes},
  author={Dai Pra, Paolo},
  journal={Applied mathematics and Optimization},
  volume={23},
  number={1},
  pages={313--329},
  year={1991},
  publisher={Springer}
}

@book{korotkov2020integrals,
  title={Integrals Related to the Error Function},
  author={Korotkov, Nikolai E and Korotkov, Alexander N},
  year={2020},
  publisher={CRC Press}
}

@book{OKSENDAL,
  title={Stochastic Differential Equations - An introduction with Applications},
  author={Oksendal, Bernt},
  year={1998},
  publisher={Springer}
}

@article{bressloff2022probabilistic,
  title={A probabilistic model of diffusion through a semi-permeable barrier},
  author={Bressloff, Paul C},
  journal={Proceedings of the Royal Society A},
  volume={478},
  number={2268},
  pages={20220615},
  year={2022},
  publisher={The Royal Society}
}

@book{GARDINER,
  title={Stochastic Methods-A Handbook for the Natural and Social Sciences},
  author={Gardiner, Crispin},
  year={2010},
  publisher={Springer}
}

@book{SCHUSS,
  title={Nonlinear Filtering and Optimal Phase Tracking},
  author={Schuss, Zeev},
  year={2012},
  publisher={Springer}
}

@article{BENES,
  title={Exact Finite Dimensional Filters for Certain Diffusions  with Nonlinear Drift},
  author={Benes, V},
  journal={Stochastics},
  volume={5},
  number={1-2},
  pages={65-92},
  year={1981},
  publisher={Taylor and Francis}
}

@article{CHANDRA,
  title={Stochastic Problems in Physics and Astronomy},
  author={Chandrasekhar, S},
  journal={Review of Modern Physics},
  volume={15},
  number={1},
  pages={1--88},
  year={1943},
  publisher={Wiley Online Library}
}

@book{Borodin,
  title={Handbook of Brownian Motion - Facts and Formulae},
  author={Borodin, Andrei  and Salminen, Paavo},
  year={2002},
  publisher={Birkhäuser}
}

@article{azzalini2005skew,
  title={The skew-normal distribution and related multivariate families},
  author={Azzalini, Adelchi},
  journal={Scandinavian journal of statistics},
  volume={32},
  number={2},
  pages={159--188},
  year={2005},
  publisher={Wiley Online Library}
}

@article{rogers1981markov,
  title={Markov functions},
  author={Rogers, Leonard CG and Pitman, JW},
  journal={The Annals of Probability},
  pages={573--582},
  year={1981},
  publisher={JSTOR}
}

@article{ito1963brownian,
  title={Brownian motions on a half line},
  author={It{\^o}, Kiyoshi and McKean, HP},
  journal={Illinois journal of mathematics},
  volume={7},
  number={2},
  pages={181--231},
  year={1963},
  publisher={Duke University Press}
}

@article{ouknine1991skew,
  title={“Skew-Brownian motion” and derived processes},
  author={Ouknine, Y},
  journal={Theory of Probability \& Its Applications},
  volume={35},
  number={1},
  pages={163--169},
  year={1991},
  publisher={SIAM}
}

@article{harrison1981skew,
  title={On skew Brownian motion},
  author={Harrison, John Michael and Shepp, Lawrence A},
  journal={The Annals of probability},
  pages={309--313},
  year={1981},
  publisher={JSTOR}
}

@article{gairat2017density,
  title={Density of skew Brownian motion and its functionals with application in finance},
  author={Gairat, Alexander and Shcherbakov, Vadim},
  journal={Mathematical Finance},
  volume={27},
  number={4},
  pages={1069--1088},
  year={2017},
  publisher={Wiley Online Library}
}

@article{chetrite2015nonequilibrium,
  title={Nonequilibrium Markov processes conditioned on large deviations},
  author={Chetrite, Rapha{\"e}l and Touchette, Hugo},
  journal={Annales Henri Poincar{\'e}},
  volume={16},
  number={9},
  pages={2005--2057},
  year={2015}
}

@article{kunita1969absolute,
  title={Absolute continuity of Markov processes and generators},
  author={Kunita, Hiroshi},
  journal={Nagoya Mathematical Journal},
  volume={36},
  pages={1--26},
  year={1969},
  publisher={Cambridge University Press}
}

@article{appuhamillage2011occupation,
  title={Occupation and local times for skew Brownian motion with applications to dispersion across an interface},
  author={Appuhamillage, Thilanka and Bokil, Vrushali and Thomann, Enrique and Waymire, Edward and Wood, Brian},
  journal={The Annals of Applied Probability},
  volume={21},
  number={1},
  pages={183--214},
  year={2011},
  publisher={Institute of Mathematical Statistics}
}

@article{decamps2004applications,
  title={Applications of $\delta$-function perturbation to the pricing of derivative securities},
  author={Decamps, Marc and De Schepper, Ann and Goovaerts, Marc},
  journal={Physica A: Statistical Mechanics and its Applications},
  volume={342},
  number={3-4},
  pages={677--692},
  year={2004},
  publisher={Elsevier}
}

@article{lejay2006constructions,
  title={On the constructions of the skew Brownian motion},
  author={Lejay, Antoine},
  journal={Probability Surveys},
  volume={3},
  pages={413--466},
  year={2006},
  publisher={Institute of Mathematical Statistics and Bernoulli Society}
}

@article{fleming1997asymptotics,
  title={Asymptotics for the principal eigenvalue and eigenfunction of a nearly first-order operator with large potential},
  author={Fleming, Wendell H and Sheu, Shuenn-Jyi},
  journal={The Annals of Probability},
  volume={25},
  number={4},
  pages={1953--1994},
  year={1997},
  publisher={Institute of Mathematical Statistics}
}

@article{doob1957conditional,
  title={Conditional Brownian motion and the boundary limits of harmonic functions},
  author={Doob, Joseph L},
  journal={Bulletin de la Soci{\'e}t{\'e} Math{\'e}matique de France},
  volume={85},
  pages={431--458},
  year={1957}
}

@article{ito1965transformation,
  title={Transformation of Markov processes by multiplicative functionals},
  author={Ito, Kiyosi and Watanabe, Shinzo},
  journal={Annales de l'institut Fourier},
  volume={15},
  number={1},
  pages={13--30},
  year={1965}
}

@article{azzone2022additive,
  title={Additive normal tempered stable processes for equity derivatives and power-law scaling},
  author={Azzone, Michele and Baviera, Roberto},
  journal={Quantitative Finance},
  volume={22},
  number={3},
  pages={501--518},
  year={2022},
  publisher={Taylor \& Francis}
}

@book{rogers2000diffusions,
  title={Diffusions, Markov processes and martingales: Volume 2, It{\^o} calculus},
  author={Rogers, L Chris G and Williams, David},
  volume={2},
  year={2000},
  publisher={Cambridge university press}
}

@article{carr2003finite,
  title={The finite moment log stable process and option pricing},
  author={Carr, Peter and Wu, Liuren},
  journal={The journal of finance},
  volume={58},
  number={2},
  pages={753--777},
  year={2003},
  publisher={Wiley Online Library}
}

@article{loperfido2008network,
  title={Network bias in air quality monitoring design},
  author={Loperfido, Nicola and Guttorp, Peter},
  journal={Environmetrics},
  volume={19},
  number={7},
  pages={661--671},
  year={2008},
  publisher={Wiley Online Library}
}

@article{harvey2010portfolio,
  title={Portfolio selection with higher moments},
  author={Harvey, Campbell R and Liechty, John C and Liechty, Merrill W and M{\"u}ller, Peter},
  journal={Quantitative Finance},
  volume={10},
  number={5},
  pages={469--485},
  year={2010},
  publisher={Taylor \& Francis}
}

@article{capitanio2003graphical,
  title={Graphical models for skew-normal variates},
  author={Capitanio, A and Azzalini, A and Stanghellini, Elena},
  journal={Scandinavian Journal of Statistics},
  volume={30},
  number={1},
  pages={129--144},
  year={2003},
  publisher={Wiley Online Library}
}

@book{azzalini2013skew,
  title={The skew-normal and related families},
  author={Azzalini, Adelchi},
  volume={3},
  year={2013},
  publisher={Cambridge University Press}
}

@article{massignan2014nonergodic,
  title={Nonergodic subdiffusion from Brownian motion in an inhomogeneous medium},
  author={Massignan, Pietro and Manzo, Carlo and Torreno-Pina, Juan A and Garc{\'\i}a-Parajo, Maria F and Lewenstein, Maciej and Lapeyre Jr, Gerald J},
  journal={Physical review letters},
  volume={112},
  number={15},
  pages={150603},
  year={2014},
  publisher={APS}
}

@article{alves2016characterization,
  title={Characterization of diffusion processes: Normal and anomalous regimes},
  author={Alves, Samuel B and de Oliveira Jr, Gilson F and de Oliveira, Luimar C and de Silans, Thierry Passerat and Chevrollier, Martine and Ori{\'a}, Marcos and Cavalcante, Hugo LD de S},
  journal={Physica A: Statistical Mechanics and its Applications},
  volume={447},
  pages={392--401},
  year={2016},
  publisher={Elsevier}
}

@book{hayashi2011econometrics,
  title={Econometrics},
  author={Hayashi, Fumio},
  year={2011},
  publisher={Princeton University Press}
}

@article{azzalini1985class,
  title={A class of distributions which includes the normal ones},
  author={Azzalini, Adelchi},
  journal={Scandinavian journal of statistics},
  pages={171--178},
  year={1985},
  publisher={JSTOR}
}

@article{heckman1979sample,
  title={Sample selection bias as a specification error},
  author={Heckman, James J},
  journal={Econometrica},
  pages={153--161},
  year={1979},
  publisher={JSTOR}
}

@article{azzalini2019sample,
  title={Sample selection models for discrete and other non-Gaussian response variables},
  author={Azzalini, Adelchi and Kim, Hyoung-Moon and Kim, Hea-Jung},
  journal={Statistical Methods \& Applications},
  volume={28},
  number={1},
  pages={27--56},
  year={2019},
  publisher={Springer}
}

@article{fleming1977exit,
  title={Exit probabilities and optimal stochastic control},
  author={Fleming, Wendell H},
  journal={Applied Mathematics and Optimization},
  volume={4},
  pages={329--346},
  year={1977},
  publisher={Springer}
}

@article{jha2023skew,
  title={A Skew-Normal Spatial Simultaneous Autoregressive Model and its Implementation},
  author={Jha, Sanjeeva Kumar and Singh, Ningthoukhongjam Vikimchandra},
  journal={Sankhya A},
  volume={85},
  number={1},
  pages={306--323},
  year={2023},
  publisher={Springer}
}

@article{marchenko2012heckman,
  title={A {H}eckman selection-t model},
  author={Marchenko, Yulia V and Genton, Marc G},
  journal={Journal of the American Statistical Association},
  volume={107},
  number={497},
  pages={304--317},
  year={2012},
  publisher={Taylor \& Francis Group}
}


\section*{Appendix A: Proof of Proposition 2}

\noindent Adopt the notations:

\begin{equation}
\label{A1}
\begin{array}{l}
\Pi^{\pm} [x,  t ] = G[1\pm E], \\\\
G:= \frac{e^{- \frac{x^{2}}{2t} }} {\sqrt{2 \pi t}} \quad{\rm and}\quad E:=  {\rm Erf }[Z(t)x]= \frac{2}{\sqrt{\pi}}\int_0^{Z(t) x} e^{-s^{2}} ds, \\\\
e = \frac{2}{\sqrt{\pi}} e^{- Z^{2}(t) x^{2} } .
\end{array}
\end{equation}

\noindent Writing for short   $Z:= Z(t)$,  we have:

\begin{equation}
\label{A2}
\begin{array}{l}
\partial_t G = \frac{1}{2} \partial_{xx}G, \quad \partial_x G = -\frac{x} {t} G \quad  \partial_x e = -2Z^{2}xe \\\\
\partial_t E = \dot{Z} x e , \quad \partial_x E = Z e, \quad  \partial_{xx}E = Z \partial_x e= - 2Z^{3} x  e,\\\\
\end{array}
\end{equation}

\noindent The pde Eq.(\ref{MATCH1}) is rewritten, (omitting the arguments):

\begin{equation}
\label{A3}
\partial_t \Pi^{\pm} = - \psi \partial_x \left[ \frac{ \pm Z e}{1\pm E}  \Pi^{\pm} \right] + \frac{1}{2} \partial_{xx} \Pi^{\pm}.
\end{equation}

\noindent Using Eqs.(\ref{A1}) and (\ref{A2}), Eq.(\ref{A3}) reads:

\begin{equation}
\label{A4}
\begin{array}{l}
 \textcolor{red}{[  \partial_t G][1\pm E]  } +[\dot{Z} x e] G = - \psi \partial_x \left[ \pm Z e G\right] +\textcolor{red}{ \frac{1}{2} [1\pm E] \partial_{xx} G } \pm Z e\partial_x G \mp (Z^{3} x e) G.
\end{array}
\end{equation}

\noindent Since $\partial_t G = \frac{1}{2} \partial_{xx} G$, Eq.(\ref{A4}) can finally  be rewritten as:

\begin{equation}
\label{A5}
\begin{array}{l}
\dot{Z} = 
  Z   \frac{ [\psi -1 ] }{t} 
  + \left[ 2\psi - 1\right] Z^{3} .
  \end{array}
\end{equation}

\noindent Introducing  $Z^{2} = \frac{\rho^{2} }{2t(1-\rho^{2})}$ into Eq.(\ref{A5}), the assertion Eq.(\ref{RHOPROP})  follows.

\rightline{$\Box$}


\section*{Appendix B: The ${\rm Erfi}(x)$ change of measure}

\noindent One may wonder whether it exists a  diffusion similar to Eq.(\ref{SDESKEW}) but defined  for  $ t \in \mathbb{R}^{+}$. To address this issue, one notes that  the imaginary error function  ${\rm Erfi}(x)$ might be used, since we have:

\begin{equation}
\label{IMAGE}
\left\{
\begin{array}{l}
{\rm Erfi}\left[ \frac{x}{\sqrt{2(T+t)}}\right] := \frac{2}{\sqrt{\pi}} \int_{0}^{\frac{x}{\sqrt{2(T+t)}}} e^{s^{2}}ds = -i {\rm Erf} \left[ \frac{ix}{\sqrt{2(T+t)}}\right],  \quad i = \sqrt{-1},\\\\
\partial_ t \left\{ {\rm Erfi}\left[ \frac{x}{\sqrt{2(T+t)}}\right] \right\}  + \frac{1}{2} \partial_{xx} \left\{ {\rm Erfi}\left[ \frac{x}{\sqrt{2(T+t)}}\right] \right\}=0.
\end{array}
\right.
\end{equation}

\noindent At this stage, let us recall the  confluent hypergeometric function identity:

\begin{equation}
\label{IDENTERFI}
{\rm Erfi}(x) = \frac{2x}{\sqrt{\pi}} M\left[\frac{1}{2}, \frac{3}{2}, x^{2} \right] := \frac{2x}{\sqrt{\pi}} \sum_{n=0}^{\infty} \frac{\Gamma(a+n) \Gamma(b) }{\Gamma(b+n) \Gamma(a) }  \frac{x^{2n}}{2n!} \geq 0 \,\,\,\, {\rm for } \,\, x  \in \mathbb{R}^{+} \
\end{equation} 

\noindent  Accordingly  for  a positive definite, scalar diffusion process $X_t\in \mathbb{R}^{+}$ the function ${\rm Erfi}\left[ \frac{x}{\sqrt{2(T+t)}}\right]$ can be used as an $h$-transform. Hence  instead of Eq.\eqref{TPDSKEW}, one may tentatively write:

\begin{equation}
\label{TENTATIVE}
\left\{
\begin{array}{l}
\tilde{p}(x_0, 0|x,t) := \frac{
e^{- \frac{(x-x_0)^{2}}{2t}} } {\sqrt{2 \pi t} }
 \left\{   \frac{1 + {\rm Erfi}\left[ \frac{x}{\sqrt{2(t+T)}}\right]  }{ {1+ \rm Erfi}\left[ \frac{x_0}{\sqrt{2T}}\right]} \right\}, \qquad x, t\in \mathbb{R}^{+}, \\\\
  \displaystyle \lim_{t\rightarrow 0 } \tilde{p}(x_0, 0|x,t) = \delta(x-x_0), \qquad T>0. 
 \end{array}
 \right.
\end{equation}

\noindent Note   that $\tilde{p}(x_0, 0|x,t)$  is well defined only for $T>0$, (since 
 $\displaystyle \lim_{T\rightarrow 0 } {\rm Erfi}\left[ \frac{x_0}{\sqrt{2T}}\right]=\infty$). The Fokker-Planck equation being a continuity equation,  $\tilde{p}(x_0, 0|x,t)$ remains normalised for all $t \in \mathbb{R}^{+}$. Accordingly  $\tilde{p}(x_0, 0|x,t)$ is a regular TPD  of an underlying   diffusion process. Recalling that we have:
 
 \begin{equation}
\label{IDENT2}
{\rm Erfi}(x) \approx e^{x} \frac{1}{2x} \left[ 1 + {\rm const} \frac{1}{x}+ ...\right]  \quad {\rm for } \quad x \rightarrow \infty, 
\end{equation}
 
 \noindent  the asymptotic behaviour of  $\tilde{p}(x_0, 0|x,t)$ is given by:
 \begin{equation}
\label{ASYB}
\tilde{p}_{T}(x_0, 0|x,t) \approx \frac{\sqrt{2(T+t)}}{x} e^{- \frac{Tx^{2}}{2t(t+T)} } \quad {\rm for } \quad x \rightarrow \infty. 
\end{equation}

\noindent  Eq.(\ref{ASYB}) shows that  $T>0$  is indeed mandatory to ensure that  $\tilde{p}(x_0, 0|x,t)$  is normalisable.

\section*{Appendix C: The moments generating function}

By direct integration, we have:

$$
\label{MOMGEN}
  \begin{array}{l}
 G_{\pm}(k,t) :=   \int_{\mathbb{R}} e^{k(x-x_0)}  \frac{e^{- \frac{(x-x_0)^{2}}{2t}} \left[ 1\pm  {\rm Erf}\left( \frac{ x}{\sqrt{2 (T-t)}}\right) \right] }{\sqrt{2\pi t }\left[1\pm {\rm Erf}\left( \frac{ x_0}{\sqrt{2T}}\right)\right]} dx=_{z:= \frac{(x-x_0)}{\sqrt{2t}} } \\\\
   \frac{1 }{\sqrt{\pi  }\left[1\pm  {\rm Erf}\left( \frac{ x_0}{\sqrt{2T}}\right)\right]}
     \int_{\mathbb{R}} e^{-z^{2} + k\sqrt{2t} z }
      \left[ 1\pm {\rm Erf}\left( \frac{\sqrt{2t} z +x_0}{\sqrt{2 (T-t)}}\right) \right]dz =_{\omega:=  z - \frac{k \sqrt{t}}{\sqrt{2}}}
      \\\\
      
             \qquad \qquad \qquad \frac{e^{\frac{k^{2} t}{2}}}{\sqrt{\pi  }\left[1\pm  {\rm Erf}\left( \frac{ x_0}{\sqrt{2T}}\right)\right]}
      \int_{\mathbb{R}} e^{-\omega^{2} }
      \left[ 1\pm {\rm Erf}\left( \frac{\sqrt{2t} \omega +x_0 - kt}{\sqrt{2 (T-t)}}\right) \right]d\omega = \\\\

      \qquad \qquad \qquad \qquad \qquad \qquad \qquad \qquad \qquad \qquad 
      \frac{\left[ e^{\frac{k^{2} t}{2}}\pm {\rm Erf}\left( \frac{ x_0 -kt}{\sqrt{2T}}\right)\right] }{\left[1\pm {\rm Erf}\left( \frac{ x_0}{\sqrt{2T}}\right)\right]},
 \end{array}
$$
 
 \noindent where the last equality follows from entry 2.7.1-6 (pp. 63) of   \cite{korotkov2020integrals}.  Note that we have $G_{\pm}(0,t) =1$ as   normalisation requires.

 \section*{Appendix D: Invariance under linear transformations}

\noindent Eq.(\ref{DRIVE2}) can be rewritten as:

$$
\label{DRIVE}
\left\{
\begin{array}{l}
dX_t =  - \lambda X_t  dt+   \partial_u  \ln\left\{ \Phi(u) \right\}_{u=U_t}dt+ dW_t \qquad X_0 =x_0,\\\\ 

dU_t =  \partial_u  \ln\left\{ \Phi(u  \right\}_{u= U_t}dt+ dW_t, \, \,\, U_0 =0 \, \,\,  {\rm for} \,\,\, t\in [0,T].
\end{array}
\right.
$$

\noindent  with $\Phi(u) = 1 +  {\rm Erf}\left[ \frac{u}{\sqrt{2(T-t)}}\right]  $ and the corresponding FPE reads:

\begin{equation}
\label{BASE1}
\partial_t P =\partial_x \left\{\left[ \lambda x - \frac{\Phi_u}{\Phi} \right]  P\right\} - \partial_u  \left\{\left[  \frac{\Phi_u}{\Phi} \right] P \right\}  + \frac{1}{2} \partial_{xx} P + \partial_{xu}P + \frac{1}{2} \partial_{u} P
\end{equation}

\noindent with the notation $\Phi_u := \partial_u \Phi_u$ Write $P = \Phi Q$ and Eq.(\ref{BASE1})  reads:

\begin{equation}
\label{BASE2}
\begin{array}{l}
Q \Phi_t + \Phi  Q_t =  \Phi\partial_x \left\{ \lambda x  Q\right\}  - \Phi_u Q_x - \partial_u\left\{\Phi_u Q \right\}  + \\\\ \qquad  \qquad  \qquad 

\frac{1}{2} \Phi Q_{xx} + \Phi_u Q_x + \Phi Q_{xu} + \frac{\Phi}{2} Q_{uu} + \frac{Q}{2}  \Phi_{uu} + \Phi_u Q_u.
\end{array}
\end{equation}

\noindent After cancelling out $\Phi_u Q_x$ and developing the term  $\partial_u\left\{\Phi_u Q \right\}$ we have:

\begin{equation}
\label{BASE3}
\begin{array}{l}
Q  \Phi_t + \Phi  Q_t=  \phi\partial_x \left\{ \lambda x  Q\right\}  - Q \partial_{uu} \phi  - \partial_u Q \partial_u \phi +  \\\\ \qquad  \qquad  \qquad 

\frac{1}{2} \Phi \partial_{xx}Q + \Phi \partial_{xu} Q +  \frac{\Phi}{2} Q_{uu} + \frac{Q}{2}  \Phi_{uu} + \Phi_u Q_u.
\end{array}
\end{equation}

\noindent Cancelling out the term $\Phi_u Q_u$, factoring by $Q$ the terms $Q  \Phi_t$, $Q \partial_{uu} \phi $, $\frac{Q}{2}  \Phi_{uu}$ and the remaining terms by  $\phi$, we get:

\begin{equation}
\label{BASE4}
\begin{array}{l}
Q \left\{  \Phi_t+ \frac{1}{2} \phi_{uu}\right\}  = 
\Phi \left\{
 - Q_t +  \partial_x \left\{ \lambda x  Q\right\} + \frac{1}{2} Q_{xx} + \partial_{xu} Q_{xu} + \frac{1}{2} Q_{uu}.
 \right\} .

\end{array}
\end{equation}

\noindent The left hand side of Eq.(\ref{BASE4}) cancels out since $\Phi$ solve the backward heat equation. In the right hand side of Eq.(\ref{BASE4}), the drift term being linear in $x$,  $Q$ is  necessarily solved by a bi-variate Gaussian  probability density and it can be calculated  $Q$  explicitly. To this aim, introduce the change of variables:

\begin{equation}
\label{DRIFTLESS}
x \mapsto ye^{\lambda t} \Rightarrow  \left\{
\begin{array}{l}
\partial_t [\cdot]\mapsto  \partial_t [\cdot]+ \lambda e^{\lambda t} \partial_y [\cdot], \\\\
\partial_x [\cdot]\mapsto e^{\lambda t} \partial_y[\cdot], \\\\
\partial_u [\cdot] \mapsto \partial_u [\cdot], \\\\
Q dx  \mapsto Qe^{\lambda t} dy

\end{array}
\right.
\end{equation}

\noindent Introducing Eq.(\ref{DRIFTLESS}) into Eq.(\ref{BASE4}), we obtain:

\begin{equation}
\label{BASE5}
\partial_t Q (y,u,t|y_0, 0,0):= \partial_t Q =  \frac{e^{2\lambda t}}{2} \phi \partial_{yy}Q + e^{\lambda t} \partial_{yu} Q + \frac{1}{2} \partial_{uu}Q
\end{equation}

\noindent According to Chandrasekhar \cite{CHANDRA}, (use Eq.(295) in  \cite{CHANDRA} with the identification $\mu_1=0$ and $\mu_2= -\lambda$ ),  we have:

\begin{equation}
\label{Q}
\left\{
\begin{array}{l}
Q =\frac{ \exp \left\{ - \frac{a(t)(y-y_0)^{2} + 2h(t)(y-y_0)u+ b(t)u^{2}}{2 \Delta(t)} \right\}}{2\pi \sqrt{\Delta(t)}}, \\\\
 a(t):= t, \qquad h(t):=  \frac{\left[ 1- e^{\lambda t}\right] }{\lambda}, \qquad  b(t):= -\frac{\left[1- e^{ 2 \lambda t}  \right]}{2\lambda}, \qquad \Delta(t): = ab-h^{2}.
\end{array}
\right.
\end{equation}

\noindent Hence we end with:

\begin{equation}
\label{P}
P (y,u,t|y_0, 0,0) = Q \phi =  Q  \left[ 1+ {\rm Erf} \left(\frac{u}{\sqrt{2(T-t)} }\right)  \right]  dy
\end{equation}

\noindent  Calculate finally  the marginal density $P_M (y,t|y_0, 0)$ using entry 2.6.1-6  of \cite{korotkov2020integrals}  to get:
\begin{equation}
\label{MARGO}
P_M (y,t|y_0, 0,0) = \int_{\mathbb{R}} P (y,u,t|y_0, 0,0) du =  \frac{e^{- \frac{ (y-y_0)^{2} } {2 b(t)}} }{\sqrt{2 \pi b(t)}} \left\{ 1 + {\rm Erf} \left[ \frac{-h(t)(y-y_0)}{\sqrt{2[b(t)^{2} (T-t)+ b(t) \Delta(t)]}} \right] \right\} .
\end{equation}

\noindent Returning to the nominal variables for $P_M (x,t|x_0, 0)$ we obtain Eq.(\ref{MARGO1}).

\end{document}